\documentclass[conference,letterpaper]{IEEEtran}
\usepackage[cmex10]{amsmath}
\usepackage{amssymb}
\usepackage{mathtools}
\usepackage{graphicx}	
\usepackage{tikz}
\usepackage{pgfplots}
\usepackage{algorithmic}
\usepackage[]{algorithm2e}
\usetikzlibrary{arrows}

\begin{document}
\title{Hypermatrix Representation of Unbalanced Power Distribution Systems}
\author{
\IEEEauthorblockN{Alejandro Garc\'es}
\IEEEauthorblockA{Department of Electric Power Systems \\
Universidad Tecnol\'ogica de Pereira\\ 
Pereira, Colombia\\ alejandro.garces@utp.edu.co}}

\maketitle
\begin{abstract}
Tensor analysis has been a widely studied in physics applications including circuit theory and electric machines.   This paper reviews some of the main features of this type of representation for unbalanced power distribution systems and discusses its implementation in modern programming languages. A tensor representation simplifies the nomenclature and calculations in three phase unbalanced grids, reduces the storage requirement and improves the computation speed.  A fixed point algorithm is proposed for the power flow solution showing the flexibility of the tensor approach. Simulations in Matlab/Octave demonstrate the advantages of this representation.  
\end{abstract}



\section{Introduction}

Load-flow is a key tool for power systems operation, planning and control. This algorithm allow to determine the state variables of the system for stationary state operation.  There is a vast literature in power flow algorithms which is possible to classify into two main groups, namely, those based on the Newton's method and those based on a fixed point iteration.  The former correspond to the classic Newton-Raphson and its modifications such as the decoupled and fast-decoupled load flow algorithms \cite{el-hawary}.  The later includes the gauss-seidel algorithm and the forward-backward sweep algorithm for radial distribution systems \cite{renato}.  All these algorithm require an efficient storage of the topological characteristics of the grid by using a matrix representation.  However, this representation can be modified in order to use Tensors instead of Matrices.  This modification can be cumbersome at the beginning but presents clear advatages in unbalanced power distribution grids in the ABC reference frame.

Tensors are a key tool for modern applications where invariant laws require to be described independently of the reference frame. Although their use in power systems go back to the seminal work of Gabriel Kron \cite{kron}, they are little used in current power systems analysis. Most of the modern power systems are represented and analyzed using matrix theory.  However, unbalanced distribution systems and microgrids are applications in which tensors are more convenient, specially in respect to its implementation in script based programming languages such as Matlab, Octave or Python.  Recent attempts to include tensor analysis in power systems are \cite{flujotensor1} and \cite{flujotensor2}.  However, these works do not fully exploit the advantages of tensors in terms of nomenclature, invariance of the reference frame, geometric interpretation and computational implementation.

Power distribution grids are different from power systems in three main aspects:  lines are usually unbalanced, they lack of transposition and the relation R/X is usually high.   Transposition is the periodic swapping of phase positions of the conductors of a transmission line.  When a line has transposition, the grid can be easily transformed into independent symmetric components.  Therefore, the first aspect (unbalance) is usually studied in power systems by the use of  the method of symmetrical components, which simplifies its analysis.  This simplification is not possible in unbalanced power distribution grids (unless we accept some approximations).  In addition, the third aspect (relation R/X) makes difficult to establish approximations as in the case of power systems, especially in the formulation of the power flow.

This paper presents a tensor formulation for unbalanced power distribution systems in a modern perspective.  Special emphasis is made on the computational implementation since modern programming languanges such as Octave/Matlab, Phyton and Julia allows direct implementation of multidimensional arrays and tensors.  In addition, the paper proposes a sparse technique for storage and operation of the proposed tensors.  A fixed point iteration in the tensor space is proposed for the load flow calculation.  This methodology allows to solve the power flow in radial and meshed distribution grids.  

The rest of the paper is organized as follows:  Section II presents a brief review of tensor algebra and multidimentional arrays including basic definitions and nomenclature.  Section III shows the model for unbalanced grids and its main characteristics (specially from the point of view of the relation between nodal voltages and currents).  Next, the proposed computational implementation is described in Section IV.  A new scheme for storing the sparse tensors is also presented in this section.  After that, Section V describes some basic simulations implemented in Octave/Matlab.  Finally, conclusions and future research directions are presented in Section IV followed by relevant references.

\section{Review of Tensor Algebra and multi-dimentional arrays}

\begin{figure}
\centering
\footnotesize
\begin{tikzpicture}[x=1mm,y=1mm]
\draw[thick] (-2.5,2.5) rectangle +(5,25);
\draw[thick] (-2.5,7.5) -- +(5,0);
\draw[thick] (-2.5,12.5) -- +(5,0);
\draw[thick] (-2.5,17.5) -- +(5,0);
\draw[thick] (-2.5,22.5) -- +(5,0);
\node at (0,25) {$a_1$};
\node at (0,20) {$a_2$};
\node at (0,15) {$a_3$};
\node at (0,10) {$a_4$};
\node at (0,5) {$a_5$};
\node at (0,0) {$a_i$};

\draw[thick] (5,2.5) rectangle +(13,25);
\draw[thick] (5,7.5) -- +(13,0);
\draw[thick] (5,12.5) -- +(13,0);
\draw[thick] (5,17.5) -- +(13,0);
\draw[thick] (5,22.5) -- +(13,0);
\draw[thick] (11.5,2.5) -- +(0,25);
\node at (8,25) {$b_{11}$};
\node at (8,20) {$b_{21}$};
\node at (8,15) {$b_{31}$};
\node at (8,10) {$b_{41}$};
\node at (8,5) {$b_{51}$};
\node at (15,25) {$b_{12}$};
\node at (15,20) {$b_{22}$};
\node at (15,15) {$b_{32}$};
\node at (15,10) {$b_{42}$};
\node at (15,5) {$b_{52}$};
\node at (11.5,0) {$b_{ij}$};

\begin{scope}[xshift=24,yshift=24]
\draw[thick,fill=white] (21,2.5) rectangle +(21,15);
\draw[thick] (21,7.5) -- +(21,0);
\draw[thick] (21,12.5) -- +(21,0);
\draw[thick] (28.5,2.5) -- +(0,15);
\draw[thick] (35.5,2.5) -- +(0,15);
\node at (25,15) {$c_{113}$};
\node at (32,15) {$c_{123}$};
\node at (39,15) {$c_{133}$};
\node at (39,10) {$c_{233}$};
\node at (39,5) {$c_{333}$};
\end{scope}

\begin{scope}[xshift=12,yshift=12]
\draw[thick,fill=white] (21,2.5) rectangle +(21,15);
\draw[thick] (21,7.5) -- +(21,0);
\draw[thick] (21,12.5) -- +(21,0);
\draw[thick] (28.5,2.5) -- +(0,15);
\draw[thick] (35.5,2.5) -- +(0,15);
\node at (25,15) {$c_{112}$};
\node at (32,15) {$c_{122}$};
\node at (39,15) {$c_{132}$};
\node at (39,10) {$c_{232}$};
\node at (39,5) {$c_{332}$};
\end{scope}

\draw[thick,fill=white] (21,2.5) rectangle +(21,15);
\draw[thick] (21,7.5) -- +(21,0);
\draw[thick] (21,12.5) -- +(21,0);
\draw[thick] (28.5,2.5) -- +(0,15);
\draw[thick] (35.5,2.5) -- +(0,15);
\node at (25,15) {$c_{111}$};
\node at (25,10) {$c_{211}$};
\node at (25,5) {$c_{311}$};
\node at (32,15) {$c_{121}$};
\node at (32,10) {$c_{221}$};
\node at (32,5) {$c_{321}$};
\node at (39,15) {$c_{131}$};
\node at (39,10) {$c_{231}$};
\node at (39,5) {$c_{331}$};
\node at (33,0) {$c_{ijk}$};
\end{tikzpicture}
\caption{Graphical representations of a (rank-1)-tensor ($a\in \mathbb{R}^{1\times 5}$), (rank-2)-tensor ($b\in \mathbb{R}^{2\times 5}$) and (rank-3)-tensor ($c\in \mathbb{R}^{3\times 3 \times 3}$)}
\label{fig:ejemplostensores}
\end{figure}
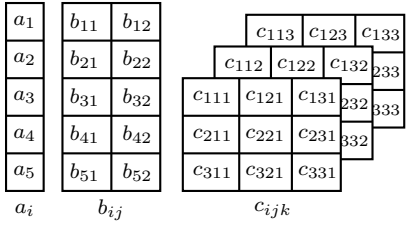

Tensors are a generalization of vectors and matrices with certain transformation properties. From computational point of view, tensors are multidimentional arrays with some useful characteristics that simplifies the analysis of complex networks \cite{Hypermatrix_algebra}.  Each tensor can be represented as a number followed by a series of subscripts.  The rank of the tensor indicates the number of subscripts to be used and the dimension of the tensor.  Thus, a rank-1 tensor is a conventional vector, a rank-2 tensor is a matrix and a rank-3 tensor is an structure similar to a rubik's cube (see Fig \ref{fig:ejemplostensores}). 

Some authors include two indices $(m,n)$ for the classification of the tensor, and allows a notation with subscripts and superscript like $x_{ij}^{kw}$.  In these representation $n$ is the number of contravariant indices (subscripts) and $m$ is the number of covariant  indices (superscpripts) \cite{librotensor1}.   This differentiation is related on how the vector changes for a given transformation. By a little abuse of notation, we will represent each tensor by a simple subscripts since in this stage, it is not relevant to differentiate between covariant and contravariant vectors.

On the other hand, tensor of the same rank can be added and multiplied by an scalar in a conventional way as

\begin{align}
    x_{ijkm} &= y_{ijkm} + z_{ijkm} \\
    x_{ijkm} &= \alpha\cdot x_{ijkm}
\end{align}

In fact, tensor can define a complete linear space in which is possible to include additional algebraic operations such as the dot product and and induced norm.   

Now, when the multiplication is between two tensors, we have expressions like this 

\begin{equation}
    y_{i} = \sum\limits_{j=1}^{N_j}\sum\limits_{k=1}^{N_k}\sum\limits_{m=1}^{N_m} x_{ijkm} \cdot y_{ijkm}
\end{equation}

This notation is greatly simplified by noticing that each repeated subscript corresponds to a sum.  In other words, the sum is assumed implicit if the index is repeated.  This results in the following clean notation

\begin{equation}
    y_{i} =  x_{ijkm} \cdot y_{ijkm}
\end{equation}

This is named the \textit{Einstein summation convention}\footnote{Some classic books include this notation only when a subscript and a superscipt is repeated.  Here, we are not differentiating between covariant and contravariant vectors allowing to use the notation even for vectors of the same type.}.

\section{Modeling Unbalanced Distribution Grids}

\subsection{Hypergraph}

\begin{figure}[tbh]
\centering
\footnotesize
\begin{tikzpicture}[x=1mm,y=1mm]
\node at (10,30)  {$\mathcal{N}=\left\{1,2,3,4 \right\}$};
\node at (18,25)  {$\mathcal{E}=\left\{1\rightarrow2,2\rightarrow3,2\rightarrow4 \right\}$};
\draw[blue,fill=blue!20] (0,0) ellipse (5 and 10);
\draw[blue,fill=blue!20] (30,0) ellipse (5 and 10);
\draw[blue,fill=blue!20] (50,30) ellipse (5 and 10);
\draw[blue,fill=blue!20] (60,0) ellipse (5 and 10);
\draw[*-] (0,5) node[left] {a} -- (30,5);
\draw[*-] (0,0) node[left] {b}-- (30,0);
\draw[*-] (0,-5) node[left] {c} -- (30,-5);

\draw[*-*] (30,5) node[above left ] {a} -- (60,5) node[right] {a};
\draw[*-*] (30,0) node[above left] {b}-- (60,0) node[right] {b};
\draw[*-*] (30,-5) node[above left] {c} -- (60,-5) node[right] {c};

\draw[-*] (30,5)  -- (50,35) node[right] {a};
\draw[-*] (30,0) -- (50,30) node[right] {b};
\draw[-*] (30,-5) -- (50,25) node[right] {c};

\node at (0,-12) {hypernode-1};
\node at (30,-12) {hypernode-2};
\node at (60,-12) {hypernode-3};
\node at (57,18) {hypernode-4};
\end{tikzpicture}
\caption{Example of an hypergraph}
\label{fig:hipergrafo}
\end{figure}
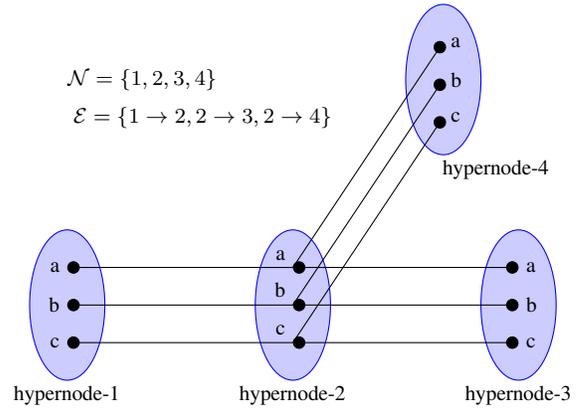

An hyper-graph is a generalization of a graph in which nodes and edges that share a similar characteristics are collected as a new set named hyper-nodes and hyper-edges as depicted in Fig \ref{fig:hipergrafo}. A three phase distribution grid can be efficiently represented by an oriented hyper-graph in which there is a set of hyper-nodes $\mathcal{N}=\left\{ 1,2,\dots,n\right\}$ each one having three different nodes $\left\{ a,b,c\right\}$. Therefore, each nodal voltage must be represented by two sub-indices as $v_{jm}$ where $j$ correspond to the phase and $m$ to the hyper-node.  In the same way, distribution lines are represented by hyper-edges $\mathcal{E}=\mathcal{N}\times\mathcal{N}$ having three different currents for each phase (or four in case to consider the neutral point).

In our representation, it is clear that each voltage and current must be stored using two sub-indices ($v_{jm}$).  This representation is also simpler for the subsequent analysis of the results.  For example, consider a three phase system with 10 nodes. Therefore, there are 30 nodal variables that require to be analyzed. In tensor representation, the voltage $v_{25}$ is clearly the voltage in the phase $b=2$ and the node $5$.  In a conventional vector representation the phases are usually stacked as $[V_{A},V_B,V_C]^T$.

\subsection{General model of the grid}
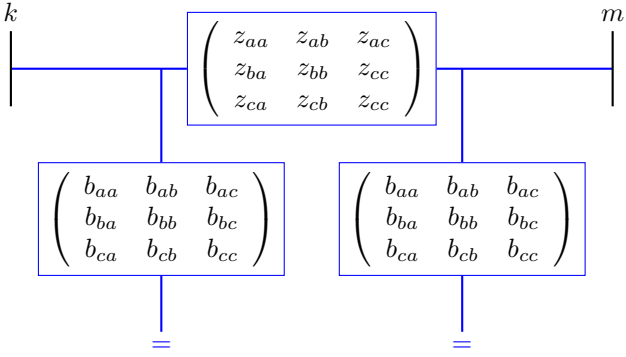
\begin{figure}
    \centering
    \begin{tikzpicture}[x=1mm,y=1mm]
    \node at (0,0) [draw=blue] (Z) {$\left( \begin{array}{ccc}z_{aa} & z_{ab} & z_{ac}\\ z_{ba} & z_{bb} & z_{cc} \\ z_{ca} & z_{cb} & z_{cc}\end{array}\right)$};
    \node at (-20,-20) [draw=blue] (B1) {$\left( \begin{array}{ccc}b_{aa} & b_{ab} & b_{ac}\\ b_{ba} & b_{bb} & b_{bc} \\ b_{ca} & b_{cb} & b_{cc}\end{array}\right)$};
    \node at (20,-20) [draw=blue] (B2) {$\left( \begin{array}{ccc}b_{aa} & b_{ab} & b_{ac}\\ b_{ba} & b_{bb} & b_{bc} \\ b_{ca} & b_{cb} & b_{cc}\end{array}\right)$};
    \draw[thick,blue] (-40,0) -- (Z) -- (40,0);
    \draw[thick,blue] (-20,0) -- (B1) -- (-20,-35) node[below] {$=$};
    \draw[thick,blue] (20,0) -- (B2) -- (20,-35) node[below] {$=$};
    \draw[thick] (-40,-5) -- +(0,10) node[above] {$k$};
    \draw[thick] (40,-5) -- +(0,10) node[above] {$m$};
    \end{tikzpicture}
    \caption{Representation of a three-phase transmission line}
    \label{fig:linea_trifasica}
\end{figure}

A feeder section in a power distribution grid can be represented by a $\pi$ model as shown in Fig \ref{fig:linea_trifasica}.  Both, series and shunt parts of the model are now $3\times 3$ complex matrices which can be non-symmetrical since, in contradistinction to the power systems, the distribution line does not have transposition.

In this case, we require four index for representing the transmission line.  For example, the term $z_{bc}$  is the mutual impedance between phase $b$ and $c$.  However, we also require to define the line section $km$; thus,  we need a (4,0)-tensor $z_{kmbc}$ in order to obtain an 
unequivocal representation.  

Using this concept, we can define a tensor for the nodal admittance ($y_{ijkm}$) which replaces the $Y_{BUS}$ in a conventional representation.  This tensor relates nodal voltages and currents as follows

\begin{equation}
i_{ik} = y_{ijkm}v_{jm}
\label{eq:tensory}
\end{equation}

where $i,j$ represent phases and $k,m$ represent nodes.  In general, currents are a covariant 2-tensor and voltages are contravariant 2-tensors, however, for our study they are simply hypermatrix. Recall here we are using the Einstein summation convention that eliminates the sum for a repeated index. This representation has interesting theoretical an practical consequences.  From the theoretical point of view it represents an invariant of the system.  This means that the expression is the same no matter the reference frame.   From the practical point of view it can be storage efficiently due to the minor symmetry of the tensor $y$, namely:

\begin{equation}
y_{ijkm} = y_{jikm} = y_{ijmk} = y_{jimk}
\label{eq:simetria}
\end{equation}

In addition, the tensor $y$ is sparse and weakly-diagonal dominant, the latear means that

\begin{equation}
\left\|y_{iikk}\right\| \geq \left\| y_{ijkm}e_{jm}\right\|
\end{equation}

where $e_{jm}$ is a vectors of ones.

\subsection{Other variables}

Tensors allows an efficient representation of the grid but also a way to consider change of basis maintaining an invariant relation.  

Let us consider a three-phase load connected to the node $m$ in which voltages and currents are given by $v_{jm}$ and $i_{jm}$ respectively.  These components can be represented in coordinates using three basis $\hat{e}_{A},\hat{e}_{B},\hat{e}_{C}$ which correspond to each phase:

\begin{equation}
    v = v_1\hat{e}_{A}+v_2\hat{e}_{B}+v_3\hat{e}_{C} 
\end{equation}

we drooped the subindex $jm$ in order to simplify the notation.  Geometrically, $v$ is just a vector represented in the space spanned by $\hat{e}_{A},\hat{e}_{B},\hat{e}_{C}$, but physically, it is a representation of the phase-to-neutral voltages.  Of course, the set of voltages in each node can be also represented in terms of the line-to-line voltages which geometrically is the same space just spanned by three new bases, namely $\hat{e}_{AB},\hat{e}_{BC},\hat{e}_{CA}$. 
Therefore, there is a set of line-to-line voltages $v'$ which can be obtained from the original phase-to-neutral voltages $v$ as follows

\begin{align}
    v_1' &= v_1'(v_1,v_2,v_3) = v_1-v_2 \\
    v_2' &= v_2'(v_1,v_2,v_3) = v_2-v_3 \\
    v_3' &= v_3'(v_1,v_2,v_3) = v_3-v_1 
\end{align}

In this case $v_j'$ represents not only the component but also the transformation given in the RHS of each equation.  Evidently, there is a linear relation between $v'$ and $v$ given by

\begin{equation}
    m_{ij} = \frac{\partial v_i'}{\partial v_j}
\end{equation}
in this case,
\begin{align}
    m_{11} &= m_{22} = m_{33} = 1 \\
    m_{12} &= m_{23} = m_{31} = -1 \\
    m_{21} &= m_{32} = m_{13} = 0
\end{align}
Therefore, line-to-line voltages can be represented as 
\begin{equation}
    v'_{jm} = m_{ij} v_{im}
\end{equation}
we can make the same for the currents, just in this case the line currents $i_{jm}$ are given in function of the currents inside a delta-connected load $i_{jm}'$:
\begin{equation}
    i_{jm} = m_{ij} i_{im}'
\end{equation}

On the other hand, the total apparent losses is by

\begin{equation}
    \bar{s}_L = \bar{v}_{ik} i_{ik} 	
\end{equation}

where $\bar{s}_L$ represents the complex conjugate of $s_L$.  This equation together with (\ref{eq:tensory}) complete the representation of the stationary state of a three-phase grid.

\section{Computational Implementation}

\subsection{Construction of the tensor $y_{ijkm}$}

The construction of the tensor $y_{ijkm}$ is  similar to the construction of the $Y_{BUS}$ matrix.  It is required only one for-loop as given in the following algorithm:

\begin{algorithm}[h]
\caption{Construction of the Tensor $y_{ijkm}$}
 \KwData{Feeder data ($z,b$), Number of nodes ($N$), Number of lines $E$}
 \KwResult{($3\times 3\times N \times N$) tensor $y_{ijkm}$}
 initialization: $y_{ijkm}\leftarrow 0+0i$\;
 
 \For{e=1 to E}{
  $\left\{ n_x,n_y\right\} \leftarrow FeederData$ \;
  $g_{(:,:)}\leftarrow Inv\left(z_{(:,:,n_x,n_y)}\right)$\;
  $y_{(:,:,n_x,n_x)}  \leftarrow y_{(:,:,n_x,n_x)} + g_{(:,:)} + b_{(:,:,n_x,n_x)}$ \;
  $y_{(:,:,n_x,n_y)}  \leftarrow y_{(:,:,n_x,n_y)} - g_{(:,:)}$ \;
  $y_{(:,:,n_y,n_x)}  \leftarrow y_{(:,:,n_y,n_x)} - g_{(:,:)}$ \;
  $y_{(:,:,n_y,n_y)}  \leftarrow y_{(:,:,n_y,n_y)} + g_{(:,:)} + b_{(:,:,n_y,n_y)}$ \;
  }
\end{algorithm}

where $g$ is a tensor such that $g_{ij} z_{ij} = \delta_{ij}$, i.e the inverse of the matrix $z$ ($\delta_{ij}$ is the Kronequer delta).  

It is important to notice that our implementation will be done in Octave/Matlab and hence the nomenclature is closely related to the actual implementation.  For example, the colon operator $(:)$ allows to retain the array shape during the assignment.  That is $y_{(:,:,n_x,n_y)}$ is a $3\times 3$ tensor which is directly assigned in the memory. Moreover, all variables are stored as a complex hyper-array.

\subsection{Storing the sparse tensor}

A $Y_{BUS}$ matrix is highly disperse. This characteristics is also inherited to the tensor $y_{ijkm}$.  In fact, this tensor has additional symmetries that enable a more efficient implementation. For this, $y_{ijkm}$ is split into two terms
\begin{equation}
	y_{ijkm} = D_{ijkm} + F_{ijkm}
\end{equation}
where $D_{ijkm}=y_{ik}\delta_{ijkm}$ being $\delta_{ijkm}$ the kronecker delta.   The storage of the tensor requires not only a reduced use of the memory but is also an efficient way to obtain back the data and to make the tensor multiplication given by (\ref{eq:tensory}). For this, let us define a structure $Y$ with six members $Y=\left\{D,F,M,C,E \right\}$. In this case, $D$ is a matrix which stores the components of ${D}_{ijkm}$ while $F$  stores the components of ${F}_{ijkm}$. Although both tensors can be stored as hyper-arrays in Octave/Matlab, the symmetry given by (\ref{eq:simetria}) allows this simpler storage.  On the other hand, $M$ is a matrix which store the non-zero positions  ($km$) of the tensor while $C$ and $E$ indicate respectively the beginning and end of each row.

\subsection{Tensor multiplication}

The proposed sparse storage allows an efficient tensor multiplication as follows:

\begin{algorithm}[h]
\caption{Tensor multiplication}
 \KwData{$Y=\left\{D,F,M,C,E \right\}, v_{jm}$}
 \KwResult{$i_{ik}$}
 initialization: $i_{ik}\leftarrow 0+0i$\;
 
 \For{i=1 to 3}{
     \For{k=1 to N}
     {
        $i_{ik} \leftarrow D_{ik}\cdot v_{ik}$ \;
     \For{$p=C_{ik}$ to $E_{ik}$}{
      $j \leftarrow M_{p1}$ \;
      $m \leftarrow M_{p2}$ \;
      $i_{ik} \leftarrow i_{ik} + F_p \cdot v_{jm}$ \; }
 }
 }
\end{algorithm}

This algorithm shows explicitly the use of the structure $Y$ and each of its data members:  $D$ and $F$ have the values themselves while $M,C$ and $E$ determine the order in which the array is  obtained.

\subsection{Load flow formulation}

\section{Numerical Results}

\begin{figure*}[bt]
\scriptsize
\centering
\begin{tikzpicture}[x=0.25mm, y = 0.25mm]
 \fill (99,-323) circle[radius=2pt] node [inner sep = 0pt](N1){}; 
 \fill (93,-284) circle[radius=2pt] node [inner sep = 0pt](N2){}; 
 \fill (99,-383) circle[radius=2pt] node [inner sep = 0pt](N3){}; 
 \fill (99,-402) circle[radius=2pt] node [inner sep = 0pt](N4){}; 
 \fill (126,-383) circle[radius=2pt] node [inner sep = 0pt](N5){}; 
 \fill (157,-384) circle[radius=2pt] node [inner sep = 0pt](N6){}; 
 \fill (123,-318) circle[radius=2pt] node [inner sep = 0pt](N7){}; 
 \fill (155,-313) circle[radius=2pt] node [inner sep = 0pt](N8){}; 
 \fill (146,-275) circle[radius=2pt] node [inner sep = 0pt](N9){}; 
 \fill (110,-272) circle[radius=2pt] node [inner sep = 0pt](N10){}; 
 \fill (64,-254) circle[radius=2pt] node [inner sep = 0pt](N11){}; 
 \fill (140,-337) circle[radius=2pt] node [inner sep = 0pt](N12){}; 
 \fill (187,-307) circle[radius=2pt] node [inner sep = 0pt](N13){}; 
 \fill (124,-249) circle[radius=2pt] node [inner sep = 0pt](N14){}; 
 \fill (207,-369) circle[radius=2pt] node [inner sep = 0pt](N15){}; 
 \fill (219,-397) circle[radius=2pt] node [inner sep = 0pt](N16){}; 
 \fill (234,-360) circle[radius=2pt] node [inner sep = 0pt](N17){}; 
 \fill (149,-184) circle[radius=2pt] node [inner sep = 0pt](N18){}; 
 \fill (111,-193) circle[radius=2pt] node [inner sep = 0pt](N19){}; 
 \fill (71,-206) circle[radius=2pt] node [inner sep = 0pt](N20){}; 
 \fill (137,-146) circle[radius=2pt] node [inner sep = 0pt](N21){}; 
 \fill (63,-165) circle[radius=2pt] node [inner sep = 0pt](N22){}; 
 \fill (123,-101) circle[radius=2pt] node [inner sep = 0pt](N23){}; 
 \fill (65,-116) circle[radius=2pt] node [inner sep = 0pt](N24){}; 
 \fill (114,-70) circle[radius=2pt] node [inner sep = 0pt](N25){}; 
 \fill (76,-78) circle[radius=2pt] node [inner sep = 0pt](N26){}; 
 \fill (31,-87) circle[radius=2pt] node [inner sep = 0pt](N27){}; 
 \fill (107,-44) circle[radius=2pt] node [inner sep = 0pt](N28){}; 
 \fill (99,-17) circle[radius=2pt] node [inner sep = 0pt](N29){}; 
 \fill (130,-14) circle[radius=2pt] node [inner sep = 0pt](N30){}; 
 \fill (64,-44) circle[radius=2pt] node [inner sep = 0pt](N31){}; 
 \fill (55,-17) circle[radius=2pt] node [inner sep = 0pt](N32){}; 
 \fill (22,-42) circle[radius=2pt] node [inner sep = 0pt](N33){}; 
 \fill (197,-337) circle[radius=2pt] node [inner sep = 0pt](N34){}; 
 \fill (204,-171) circle[radius=2pt] node [inner sep = 0pt](N35){}; 
 \fill (260,-188) circle[radius=2pt] node [inner sep = 0pt](N36){}; 
 \fill (207,-202) circle[radius=2pt] node [inner sep = 0pt](N37){}; 
 \fill (294,-176) circle[radius=2pt] node [inner sep = 0pt](N38){}; 
 \fill (322,-168) circle[radius=2pt] node [inner sep = 0pt](N39){}; 
 \fill (195,-148) circle[radius=2pt] node [inner sep = 0pt](N40){}; 
 \fill (244,-133) circle[radius=2pt] node [inner sep = 0pt](N41){}; 
 \fill (190,-124) circle[radius=2pt] node [inner sep = 0pt](N42){}; 
 \fill (245,-104) circle[radius=2pt] node [inner sep = 0pt](N43){}; 
 \fill (181,-98) circle[radius=2pt] node [inner sep = 0pt](N44){}; 
 \fill (218,-88) circle[radius=2pt] node [inner sep = 0pt](N45){}; 
 \fill (258,-75) circle[radius=2pt] node [inner sep = 0pt](N46){}; 
 \fill (170,-66) circle[radius=2pt] node [inner sep = 0pt](N47){}; 
 \fill (141,-74) circle[radius=2pt] node [inner sep = 0pt](N48){}; 
 \fill (209,-54) circle[radius=2pt] node [inner sep = 0pt](N49){}; 
 \fill (245,-44) circle[radius=2pt] node [inner sep = 0pt](N50){}; 
 \fill (276,-36) circle[radius=2pt] node [inner sep = 0pt](N51){}; 
 \fill (268,-293) circle[radius=2pt] node [inner sep = 0pt](N52){}; 
 \fill (298,-288) circle[radius=2pt] node [inner sep = 0pt](N53){}; 
 \fill (323,-282) circle[radius=2pt] node [inner sep = 0pt](N54){}; 
 \fill (350,-281) circle[radius=2pt] node [inner sep = 0pt](N55){}; 
 \fill (377,-276) circle[radius=2pt] node [inner sep = 0pt](N56){}; 
 \fill (313,-238) circle[radius=2pt] node [inner sep = 0pt](N57){}; 
 \fill (281,-247) circle[radius=2pt] node [inner sep = 0pt](N58){}; 
 \fill (245,-255) circle[radius=2pt] node [inner sep = 0pt](N59){}; 
 \fill (381,-224) circle[radius=2pt] node [inner sep = 0pt](N60){}; 
 \fill (394,-256) circle[radius=2pt] node [inner sep = 0pt](N61){}; 
 \fill (380,-154) circle[radius=2pt] node [inner sep = 0pt](N62){}; 
 \fill (369,-119) circle[radius=2pt] node [inner sep = 0pt](N63){}; 
 \fill (358,-90) circle[radius=2pt] node [inner sep = 0pt](N64){}; 
 \fill (318,-107) circle[radius=2pt] node [inner sep = 0pt](N65){}; 
 \fill (326,-133) circle[radius=2pt] node [inner sep = 0pt](N66){}; 
 \fill (462,-213) circle[radius=2pt] node [inner sep = 0pt](N67){}; 
 \fill (485,-200) circle[radius=2pt] node [inner sep = 0pt](N68){}; 
 \fill (510,-186) circle[radius=2pt] node [inner sep = 0pt](N69){}; 
 \fill (534,-171) circle[radius=2pt] node [inner sep = 0pt](N70){}; 
 \fill (560,-152) circle[radius=2pt] node [inner sep = 0pt](N71){}; 
 \fill (475,-250) circle[radius=2pt] node [inner sep = 0pt](N72){}; 
 \fill (503,-236) circle[radius=2pt] node [inner sep = 0pt](N73){}; 
 \fill (527,-220) circle[radius=2pt] node [inner sep = 0pt](N74){}; 
 \fill (552,-206) circle[radius=2pt] node [inner sep = 0pt](N75){}; 
 \fill (487,-293) circle[radius=2pt] node [inner sep = 0pt](N76){}; 
 \fill (505,-285) circle[radius=2pt] node [inner sep = 0pt](N77){}; 
 \fill (522,-277) circle[radius=2pt] node [inner sep = 0pt](N78){}; 
 \fill (540,-266) circle[radius=2pt] node [inner sep = 0pt](N79){}; 
 \fill (524,-308) circle[radius=2pt] node [inner sep = 0pt](N80){}; 
 \fill (527,-343) circle[radius=2pt] node [inner sep = 0pt](N81){}; 
 \fill (528,-383) circle[radius=2pt] node [inner sep = 0pt](N82){}; 
 \fill (560,-383) circle[radius=2pt] node [inner sep = 0pt](N83){}; 
 \fill (563,-318) circle[radius=2pt] node [inner sep = 0pt](N84){}; 
 \fill (561,-248) circle[radius=2pt] node [inner sep = 0pt](N85){}; 
 \fill (492,-361) circle[radius=2pt] node [inner sep = 0pt](N86){}; 
 \fill (437,-365) circle[radius=2pt] node [inner sep = 0pt](N87){}; 
 \fill (429,-324) circle[radius=2pt] node [inner sep = 0pt](N88){}; 
 \fill (395,-369) circle[radius=2pt] node [inner sep = 0pt](N89){}; 
 \fill (391,-332) circle[radius=2pt] node [inner sep = 0pt](N90){}; 
 \fill (352,-375) circle[radius=2pt] node [inner sep = 0pt](N91){}; 
 \fill (349,-341) circle[radius=2pt] node [inner sep = 0pt](N92){}; 
 \fill (312,-376) circle[radius=2pt] node [inner sep = 0pt](N93){}; 
 \fill (304,-326) circle[radius=2pt] node [inner sep = 0pt](N94){}; 
 \fill (280,-381) circle[radius=2pt] node [inner sep = 0pt](N95){}; 
 \fill (274,-333) circle[radius=2pt] node [inner sep = 0pt](N96){}; 
 \fill (447,-174) circle[radius=2pt] node [inner sep = 0pt](N97){}; 
 \fill (467,-160) circle[radius=2pt] node [inner sep = 0pt](N98){}; 
 \fill (487,-148) circle[radius=2pt] node [inner sep = 0pt](N99){}; 
 \fill (522,-122) circle[radius=2pt] node [inner sep = 0pt](N100){}; 
 \fill (434,-136) circle[radius=2pt] node [inner sep = 0pt](N101){}; 
 \fill (459,-119) circle[radius=2pt] node [inner sep = 0pt](N102){}; 
 \fill (489,-101) circle[radius=2pt] node [inner sep = 0pt](N103){}; 
 \fill (516,-82) circle[radius=2pt] node [inner sep = 0pt](N104){}; 
 \fill (424,-107) circle[radius=2pt] node [inner sep = 0pt](N105){}; 
 \fill (451,-89) circle[radius=2pt] node [inner sep = 0pt](N106){}; 
 \fill (492,-66) circle[radius=2pt] node [inner sep = 0pt](N107){}; 
 \fill (413,-76) circle[radius=2pt] node [inner sep = 0pt](N108){}; 
 \fill (445,-61) circle[radius=2pt] node [inner sep = 0pt](N109){}; 
 \fill (476,-37) circle[radius=2pt] node [inner sep = 0pt](N110){}; 
 \fill (429,-38) circle[radius=2pt] node [inner sep = 0pt](N111){}; 
 \fill (505,-38) circle[radius=2pt] node [inner sep = 0pt](N112){}; 
 \fill (540,-37) circle[radius=2pt] node [inner sep = 0pt](N113){}; 
 \fill (572,-37) circle[radius=2pt] node [inner sep = 0pt](N114){}; 
 \fill (70,-323) circle[radius=2pt] node [inner sep = 0pt](N116){}; 
 \fill (30,-323) circle[radius=2pt] node [inner sep = 0pt](N117){}; 
 \fill (420,-220) circle[radius=2pt] node [inner sep = 0pt](N120){}; 
 \fill (180,-175) circle[radius=2pt] node [inner sep = 0pt](N115){}; 
 \fill (220,-300) circle[radius=2pt] node [inner sep = 0pt](N119){}; 
 \fill (420,-256) circle[radius=2pt] node [inner sep = 0pt](N129){}; 
 \fill (441,-155) circle[radius=2pt] node [inner sep = 0pt](N122){}; 
 \fill (310,-36) circle[radius=2pt] node [inner sep = 0pt](N118){}; 
 \fill (150,-14) circle[radius=2pt] node [inner sep = 0pt](N123){}; 
 \fill (360,-36) circle[radius=2pt] node [inner sep = 0pt](N125){}; 
 \fill (550,-122) circle[radius=2pt] node [inner sep = 0pt](N127){}; 
 \node at (N1)[yshift=-6,xshift=-6] {1}; 
 \node at (N2)[yshift=6] {2}; 
 \node at (N3)[xshift=-6] {3}; 
 \node at (N4)[yshift=-6] {4}; 
 \node at (N5)[yshift=6] {5}; 
 \node at (N6)[yshift=6] {6}; 
 \node at (N7)[yshift=6] {7}; 
 \node at (N8)[yshift=-6] {8}; 
 \node at (N9)[yshift=6] {9}; 
 \node at (N10)[yshift=-6] {10}; 
 \node at (N11)[yshift=6] {11}; 
 \node at (N12)[yshift=-6] {12}; 
 \node at (N13)[yshift=6,xshift=6] {13}; 
 \node at (N14)[yshift=6] {14}; 
 \node at (N15)[xshift=-6] {15}; 
 \node at (N16)[yshift=-6] {16}; 
 \node at (N17)[yshift=6] {17}; 
 \node at (N18)[yshift=6,xshift=4] {18}; 
 \node at (N19)[yshift=6] {19}; 
 \node at (N20)[yshift=-6] {20}; 
 \node at (N21)[xshift=6] {21}; 
 \node at (N22)[yshift=6] {22}; 
 \node at (N23)[xshift=6] {23}; 
 \node at (N24)[yshift=6] {24}; 
 \node at (N25)[xshift=6] {25}; 
 \node at (N26)[yshift=-6] {26}; 
 \node at (N27)[yshift=-6] {27}; 
 \node at (N28)[xshift=6] {28}; 
 \node at (N29)[yshift=6] {29}; 
 \node at (N30)[yshift=-6] {30}; 
 \node at (N31)[xshift=6] {31}; 
 \node at (N32)[yshift=6] {32}; 
 \node at (N33)[yshift=6] {33}; 
 \node at (N34)[xshift=6] {34}; 
 \node at (N35)[xshift=6,yshift=5] {35}; 
 \node at (N36)[yshift=6] {36}; 
 \node at (N37)[yshift=6] {37}; 
 \node at (N38)[yshift=6] {38}; 
 \node at (N39)[yshift=-6] {39}; 
 \node at (N40)[xshift=-6] {40}; 
 \node at (N41)[yshift=6] {41}; 
 \node at (N42)[xshift=-6] {42}; 
 \node at (N43)[yshift=6] {43}; 
 \node at (N44)[xshift=-6] {44}; 
 \node at (N45)[yshift=6] {45}; 
 \node at (N46)[yshift=6] {46}; 
 \node at (N47)[yshift=6] {47}; 
 \node at (N48)[yshift=6] {48}; 
 \node at (N49)[yshift=6] {49}; 
 \node at (N50)[yshift=-6] {50}; 
 \node at (N51)[yshift=6] {51}; 
 \node at (N52)[yshift=6] {52}; 
 \node at (N53)[yshift=6] {53}; 
 \node at (N54)[yshift=-6] {54}; 
 \node at (N55)[yshift=6] {55}; 
 \node at (N56)[yshift=6] {56}; 
 \node at (N57)[yshift=6] {57}; 
 \node at (N58)[yshift=6] {58}; 
 \node at (N59)[yshift=6] {59}; 
 \node at (N60)[yshift=6,xshift=6] {60}; 
 \node at (N61)[yshift=-6] {61}; 
 \node at (N62)[xshift=6] {62}; 
 \node at (N63)[xshift=6] {63}; 
 \node at (N64)[yshift=6] {64}; 
 \node at (N65)[yshift=6] {65}; 
 \node at (N66)[yshift=-6] {66}; 
 \node at (N67)[yshift=6,xshift=-8] {67}; 
 \node at (N68)[yshift=6] {68}; 
 \node at (N69)[yshift=6] {69}; 
 \node at (N70)[yshift=6] {70}; 
 \node at (N71)[yshift=6] {71}; 
 \node at (N72)[xshift=-6] {72}; 
 \node at (N73)[yshift=6] {73}; 
 \node at (N74)[yshift=6] {74}; 
 \node at (N75)[yshift=6] {75}; 
 \node at (N76)[xshift=-6] {76}; 
 \node at (N77)[yshift=6] {77}; 
 \node at (N78)[yshift=6] {78}; 
 \node at (N79)[yshift=6] {79}; 
 \node at (N80)[xshift=6] {80}; 
 \node at (N81)[xshift=-6] {81}; 
 \node at (N82)[xshift=-6] {82}; 
 \node at (N83)[yshift=6] {83}; 
 \node at (N84)[yshift=-6] {84}; 
 \node at (N85)[yshift=6] {85}; 
 \node at (N86)[yshift=6, xshift=-7] {86}; 
 \node at (N87)[yshift=-6] {87}; 
 \node at (N88)[yshift=6] {88}; 
 \node at (N89)[yshift=-6] {89}; 
 \node at (N90)[yshift=6] {90}; 
 \node at (N91)[yshift=-6] {91}; 
 \node at (N92)[yshift=6] {92}; 
 \node at (N93)[yshift=-6] {93}; 
 \node at (N94)[yshift=6] {94}; 
 \node at (N95)[yshift=-6] {95}; 
 \node at (N96)[yshift=6] {96}; 
 \node at (N97)[xshift=-7] {97}; 
 \node at (N98)[yshift=6] {98}; 
 \node at (N99)[yshift=6] {99}; 
 \node at (N100)[yshift=6] {100}; 
 \node at (N101)[xshift=-9] {101}; 
 \node at (N102)[yshift=6] {102}; 
 \node at (N103)[yshift=6] {103}; 
 \node at (N104)[yshift=6] {104}; 
 \node at (N105)[xshift=-9] {105}; 
 \node at (N106)[yshift=6] {106}; 
 \node at (N107)[yshift=6] {107}; 
 \node at (N108)[yshift=6] {108}; 
 \node at (N109)[yshift=6] {109}; 
 \node at (N110)[yshift=6] {110}; 
 \node at (N111)[yshift=6] {111}; 
 \node at (N112)[yshift=6] {112}; 
 \node at (N113)[yshift=6] {113}; 
 \node at (N114)[yshift=6] {114}; 
 \node at (N116)[yshift=6] {149}; 
 \node at (N117)[yshift=6] {150}; 
 \node at (N120)[yshift=6] {160}; 
 \node at (N115)[yshift=6] {135}; 
 \node at (N119)[yshift=6] {152}; 
 \node at (N129)[yshift=6] {610}; 
 \node at (N122)[xshift=-9] {197}; 
 \node at (N118)[yshift=6] {151}; 
 \node at (N123)[yshift=6] {250}; 
 \node at (N125)[yshift=6] {300}; 
 \node at (N127)[yshift=6] {450};

 \draw[thick,black] (N116) -- (N1); 
 \draw[thick,gray] (N1) -- (N2); 
 \draw[thick,gray] (N1) -- (N3); 
 \draw[thick,gray] (N3) -- (N4); 
 \draw[thick,gray] (N3) -- (N5); 
 \draw[thick,gray] (N5) -- (N6); 
 \draw[thick,black] (N1) -- (N7); 
 \draw[thick,black] (N7) -- (N8); 
 \draw[thick,gray] (N8) -- (N9); 
 \draw[thick,gray] (N14) -- (N10); 
 \draw[thick,gray] (N14) -- (N11); 
 \draw[thick,gray] (N8) -- (N12); 
 \draw[thick,black] (N8) -- (N13); 
 \draw[thick,gray] (N9) -- (N14); 
 \draw[thick,gray] (N34) -- (N15); 
 \draw[thick,gray] (N15) -- (N16); 
 \draw[thick,gray] (N15) -- (N17); 
 \draw[thick,black] (N13) -- (N18); 
 \draw[thick,gray] (N18) -- (N19); 
 \draw[thick,gray] (N19) -- (N20); 
 \draw[thick,black] (N18) -- (N21); 
 \draw[thick,gray] (N21) -- (N22); 
 \draw[thick,black] (N21) -- (N23); 
 \draw[thick,gray] (N23) -- (N24); 
 \draw[thick,black] (N23) -- (N25); 
 \draw[thick,gray] (N25) -- (N26); 
 \draw[thick,gray] (N26) -- (N27); 
 \draw[thick,black] (N25) -- (N28); 
 \draw[thick,black] (N28) -- (N29); 
 \draw[thick,black] (N29) -- (N30); 
 \draw[thick,gray] (N26) -- (N31); 
 \draw[thick,gray] (N31) -- (N32); 
 \draw[thick,gray] (N27) -- (N33); 
 \draw[thick,gray] (N13) -- (N34); 
 \draw[thick,black] (N115) -- (N35); 
 \draw[thick,gray] (N35) -- (N36); 
 \draw[thick,gray] (N36) -- (N37); 
 \draw[thick,gray] (N36) -- (N38); 
 \draw[thick,gray] (N38) -- (N39); 
 \draw[thick,black] (N35) -- (N40); 
 \draw[thick,gray] (N40) -- (N41); 
 \draw[thick,black] (N40) -- (N42); 
 \draw[thick,gray] (N42) -- (N43); 
 \draw[thick,black] (N42) -- (N44); 
 \draw[thick,gray] (N44) -- (N45); 
 \draw[thick,gray] (N45) -- (N46); 
 \draw[thick,black] (N44) -- (N47); 
 \draw[thick,black] (N47) -- (N48); 
 \draw[thick,black] (N47) -- (N49); 
 \draw[thick,black] (N49) -- (N50); 
 \draw[thick,black] (N50) -- (N51); 
 \draw[thick,black] (N119) -- (N52); 
 \draw[thick,black] (N52) -- (N53); 
 \draw[thick,black] (N53) -- (N54); 
 \draw[thick,black] (N54) -- (N55); 
 \draw[thick,black] (N55) -- (N56); 
 \draw[thick,black] (N54) -- (N57); 
 \draw[thick,gray] (N57) -- (N58); 
 \draw[thick,gray] (N58) -- (N59); 
 \draw[thick,black] (N57) -- (N60); 
 \draw[thick,black] (N60) -- (N61); 
 \draw[thick,black] (N60) -- (N62); 
 \draw[thick,black] (N62) -- (N63); 
 \draw[thick,black] (N63) -- (N64); 
 \draw[thick,black] (N64) -- (N65); 
 \draw[thick,black] (N65) -- (N66); 
 \draw[thick,black] (N120) -- (N67); 
 \draw[thick,gray] (N67) -- (N68); 
 \draw[thick,gray] (N68) -- (N69); 
 \draw[thick,gray] (N69) -- (N70); 
 \draw[thick,gray] (N70) -- (N71); 
 \draw[thick,black] (N67) -- (N72); 
 \draw[thick,gray] (N72) -- (N73); 
 \draw[thick,gray] (N73) -- (N74); 
 \draw[thick,gray] (N74) -- (N75); 
 \draw[thick,black] (N72) -- (N76); 
 \draw[thick,black] (N76) -- (N77); 
 \draw[thick,black] (N77) -- (N78); 
 \draw[thick,black] (N78) -- (N79); 
 \draw[thick,black] (N78) -- (N80); 
 \draw[thick,black] (N80) -- (N81); 
 \draw[thick,black] (N81) -- (N82); 
 \draw[thick,black] (N82) -- (N83); 
 \draw[thick,gray] (N81) -- (N84); 
 \draw[thick,gray] (N84) -- (N85); 
 \draw[thick,black] (N76) -- (N86); 
 \draw[thick,black] (N86) -- (N87); 
 \draw[thick,gray] (N87) -- (N88); 
 \draw[thick,black] (N87) -- (N89); 
 \draw[thick,gray] (N89) -- (N90); 
 \draw[thick,black] (N89) -- (N91); 
 \draw[thick,gray] (N91) -- (N92); 
 \draw[thick,black] (N91) -- (N93); 
 \draw[thick,gray] (N93) -- (N94); 
 \draw[thick,black] (N93) -- (N95); 
 \draw[thick,gray] (N95) -- (N96); 
 \draw[thick,black] (N67) -- (N97); 
 \draw[thick,black] (N97) -- (N98); 
 \draw[thick,black] (N98) -- (N99); 
 \draw[thick,black] (N99) -- (N100); 
 \draw[thick,black] (N122) -- (N101); 
 \draw[thick,gray] (N101) -- (N102); 
 \draw[thick,gray] (N102) -- (N103); 
 \draw[thick,gray] (N103) -- (N104); 
 \draw[thick,black] (N101) -- (N105); 
 \draw[thick,gray] (N105) -- (N106); 
 \draw[thick,gray] (N106) -- (N107); 
 \draw[thick,black] (N105) -- (N108); 
 \draw[thick,gray] (N108) -- (N109); 
 \draw[thick,gray] (N109) -- (N110); 
 \draw[thick,gray] (N110) -- (N111); 
 \draw[thick,gray] (N110) -- (N112); 
 \draw[thick,gray] (N112) -- (N113); 
 \draw[thick,gray] (N113) -- (N114); 
 \draw[thick,black] (N51) -- (N118); 
 \draw[thick,black] (N30) -- (N123); 
 \draw[thick,black] (N108) -- (N125); 
 \draw[thick,black] (N100) -- (N127); 
 \draw[thick,gray] (N61) -- (N129); 
 \draw[thick,black] (N18) -- (N115); 
 \draw[thick,black] (N117) -- (N116); 
 \draw[thick,black] (N13) -- (N119); 
 \draw[thick,black] (N60) -- (N120); 
 \draw[thick,black] (N97) -- (N122); 
\node at (N117) [yshift=-8,blue!50!green]  {Slack};
\end{tikzpicture}
\caption{IEEE 123 nodes test system without voltage regulators}
\label{fig:ieee123}
\end{figure*}
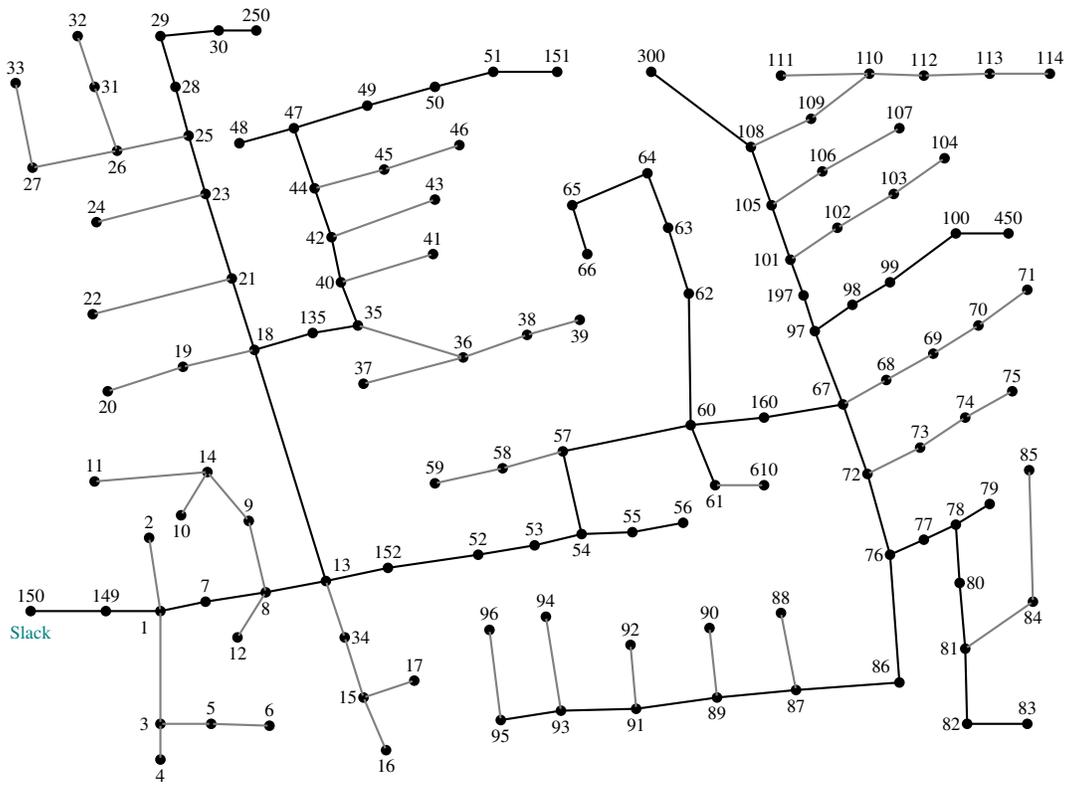

Numerical simulations where performed in the IEEE 123 test distribution system \cite{testfeeders} depicted in Fig \ref{fig:ieee123}. The system was implemented in Octave/Matlab and executed in a computer with $i7-6700 CPU/3.4 GHz$. The code in Octave/Matlab is available in \cite{Matlab_exchange}.  
A simple calculation for calculation of voltages and currents was performed using a three phase $y_{BUS}$ and the proposed tensor representation with sparse implementation. The size of each array and the elapsed time for obtaining currents from the voltages are given in Table \ref{tab:resultados}.   There is a reduction of more than $95\%$ in the execution time by using the proposed tensor representation with sparse storage.  This is an important reduction in many practical problems in which this sub-rutine is executed many times.

\begin{table}[tb]
\centering
\caption{Numerical results for different implementations of the grid model}
\label{tab:resultados}
\begin{tabular}{|c|c|c|}
\hline
Parameter & $Y_{BUS}$ (matrix) &
$y_{ijkm}$ \\
\hline\hline
Elapsed time & 0.151272 &   0.005895 \\\hline
Array's size & $357\times 357$ & $3\times 3 \times 119\times 119$  \\
\hline     
\end{tabular}
\end{table}

\begin{table}[tb]
    \centering
    \caption{Parameters of the structure Y}
    \label{tab:parametros_estructura}
    \begin{tabular}{|c|c|}
    \hline
    Member of the structure & size \\
    \hline\hline
   $M$ & $1586\times 2$ \\
   $D$ & $3\times 119$ \\
   $F$ & $1\times 1586$ \\
   $C$ & $3\times 119$ \\
   $E$ & $3\times 119$ \\
\hline
    \end{tabular}
\end{table}

Sparse representation of the tensor required only 5829 memory positions while the conventional $Y_{BUS}$ representation requires 127449 memory positions.  This means, the proposed sparse representation required only $5\%$ of the memory.

\section{Conclusions}

This paper presented a tensor representation of a three-phase unbalanced distribution systems using tensors.   Some concepts of tensor algebra were reviewed and studied in the context of power distribution grids.    

The $Y_{BUS}$ matrix was represented by a $3\times\ 3\times N \times N-rank$ tensor that shows additional symmetries that are hidden in the matrix representation.   This can be used from the theoretical but practical point of view.   From the theoretical point of view, symmetries are desirable feature that allows simplifications whereas from the practical point of view, they allows a more efficient storage in digital computers.

A sparse storage of the tensor $y_{ijkm}$ was also proposed. This required only $2\%$ of the memory required  for a full matrix representation.  In addition, calculation time was only $5\%$ of the calculation for the three-phase matrix representation. This is a big reduction of the calculation time. An additional feature of the tensor representation is a 'clean code' when it is implemented in Matlab/Octave.  This may seem like a a minor advantage, but a 'clean clode'   is important in large and collaborative implementation as is usual in moder research proyects. 

Although application of tensor algebra applied to power systems goes back to the work of Gabriel Kron, this theory is little used in current power system analysis. It is the author belief that some of these concepts should be evaluated taking into account the use of programming languages such as Octave/Matlab and Python.  These languages allows to manipulate easily hyper-arrays of several complex variables.    More research is required in this direction.

\bibliographystyle{IEEEtran}
\bibliography{bibliografia}

\end{document}